\documentclass[12pt]{article}
\usepackage{latexsym,amssymb,amsmath}
\begin{document}

\baselineskip=18pt

\newcommand{\la}{\langle}
\newcommand{\ra}{\rangle}
\newcommand{\psp}{\vspace{0.4cm}}
\newcommand{\pse}{\vspace{0.2cm}}
\newcommand{\ptl}{\partial}
\newcommand{\dlt}{\delta}
\newcommand{\sgm}{\sigma}
\newcommand{\al}{\alpha}
\newcommand{\be}{\beta}
\newcommand{\G}{\Gamma}
\newcommand{\gm}{\gamma}
\newcommand{\vs}{\varsigma}
\newcommand{\Lmd}{\Lambda}
\newcommand{\lmd}{\lambda}
\newcommand{\td}{\tilde}
\newcommand{\vf}{\varphi}
\newcommand{\yt}{Y^{\nu}}
\newcommand{\wt}{\mbox{wt}\:}
\newcommand{\rd}{\mbox{Res}}
\newcommand{\ad}{\mbox{ad}}
\newcommand{\stl}{\stackrel}
\newcommand{\ol}{\overline}
\newcommand{\ul}{\underline}
\newcommand{\es}{\epsilon}
\newcommand{\dmd}{\diamond}
\newcommand{\clt}{\clubsuit}
\newcommand{\vt}{\vartheta}
\newcommand{\ves}{\varepsilon}
\newcommand{\dg}{\dagger}
\newcommand{\tr}{\mbox{Tr}}
\newcommand{\ga}{{\cal G}({\cal A})}
\newcommand{\hga}{\hat{\cal G}({\cal A})}
\newcommand{\Edo}{\mbox{End}\:}
\newcommand{\for}{\mbox{for}}
\newcommand{\kn}{\mbox{ker}}
\newcommand{\Dlt}{\Delta}
\newcommand{\rad}{\mbox{Rad}}
\newcommand{\rta}{\rightarrow}

\begin{center}{\LARGE \bf Simple Conformal Algebras Generated}\end{center}
\begin{center}{\LARGE \bf by Jordan Algebras}\end{center}

\begin{center}{\large Xiaoping Xu}\end{center}

\begin{center}{(A talk presented in Oberwolfach Conference on Jordan Algebras, 13 -- 19 August 2000)}\end{center}

\section{Background and Motivation}

We start with an example of affine Kac-Moody algebras and the Virasoro algebra. In this talk,
$\Bbb{F}$ will be a field with characteristic 0, and all the vector spaces are assumed over $\Bbb{F}$. Denote by $\Bbb{Z}$ the ring of integers and by $\Bbb{N}$ the set of nonnegative integers. 

Let $2\leq n\in\Bbb{N}$. Set
$$sl(n,\Bbb{F})=\{A\in M_{n\times n}(\Bbb{F})\mid \mbox{tr}\:A=0\},\eqno(1.1)$$
$$\la A,B\ra=\mbox{tr}\:AB\qquad\for\;\;A,B\in sl(n,\Bbb{F}),\eqno(1.2)$$
where $M_{n\times n}(\Bbb{F})$ is the algebra of $n\times n$ matrices. Then $sl(n,\Bbb{F})$ forms a simple Lie algebra with the Lie bracket
$$[A,B]=AB-BA\qquad\for\;\;A,B\in sl(n,\Bbb{F}),\eqno(1.3)$$
and $\la\cdot,\cdot\ra$ forms a symmetric invariant bilinear form, that is,
$$\la [A,B],C\ra=\la A,[B,C]\ra \qquad\for\;\;A,B\in sl(n,\Bbb{F}).\eqno(1.4)$$

Let $t$ be an indeterminant and set
$$\widehat{sl}(n,\Bbb{F})=sl(n,\Bbb{F})\otimes_{\Bbb{F}}\Bbb{F}[t,t^{-1}]\oplus\Bbb{F}\kappa\oplus \Bbb{F}d,\eqno(1.5)$$
where $\kappa$ and $d$ are symbols serving as base elements. We define a Lie bracket on $\widehat{sl}(n,\Bbb{F})$ by
$$[u\otimes t^l,v\otimes t^j]=[u,v]\otimes t^{l+j}+l\la u,v\ra\dlt_{l+j,0}\kappa,\;\;[d,u\otimes t^l]=l u\otimes t^l\eqno(1.6)$$
for $u,v\in sl(n,\Bbb{F}),\;l,j\in\Bbb{Z}$ and 
$$[\widehat{sl}(n,\Bbb{F}),\kappa]=0.\eqno(1.7)$$
The Lie algebra $(\widehat{sl}(n,\Bbb{F}),[\cdot,\cdot])$ is  an {\it affine Kac-Moody algebra}.

Define the {\it generating function}
$$u(z)=\sum_{j\in\Bbb{Z}}(u\otimes t^j)z^{-n-1}\qquad\for\;\;u\in sl(n,\Bbb{F})
,\eqno(1.8)$$
where $z$ is a formal variable. Then
$$[u(z_1),v(z_2)]=z_2^{-1}\dlt\left({z_1\over z_2}\right)[u,v](z_2)-z_2^{-1}\ptl_{z_1}\dlt\left({z_1\over z_2}\right)\la u,v\ra \kappa\eqno(1.9)$$
for $u,v\in sl(n,\Bbb{F})$, where
$$\dlt\left({z_1\over z_2}\right)=\sum_{j\in\Bbb{Z}}\left({z_1\over z_2}\right)^j.\eqno(1.10)$$

The above formula (1.9) was mathematically first used by Lepowsky and Wilson [LW] in order to study explicit integral irreducible representations of affine Kac-Moody algebras. Acting on ``highest weight'' modules, $u(z)$ for $u\in sl(n,\Bbb{F})$ are ``current fields'' in two-dimensional quantum field theory.

The {\it Virasoro algebra} ${\cal V}$ is a Lie algebra with a basis $\{L(j),\kappa\mid j\in\Bbb{Z}\}$ and the Lie bracket
$$[L(j),L(l)]=(j-l)L(j+l)+{j^3-j\over 12}\dlt_{m+n,0}\kappa,\;\;[\kappa,{\cal V}]=0\eqno(1.11)$$
for $j,l\in\Bbb{Z}$. Set
$$L(z)=\sum_{j\in\Bbb{Z}}L(j)z^{-j-2}.\eqno(1.12)$$
Then
$$[L(z_1),L(z_2)]=z_2^{-1}\dlt\left({z_1\over z_2}\right)\ptl_{z_2}L(z_2)-2z_2^{-1}\ptl_{z_1}\dlt\left({z_1\over z_2}\right)L(z_2)-{\kappa\over 12}z_2^{-1}\ptl_{z_1}^3\dlt\left({z_1\over z_2}\right)\kappa.\eqno(1.13)$$
Acting on ``highest weight'' modules, $L(z)$ is the ``energy-momentum tensor'' in two-dimensional quantum field theory.

Let us focus only on the parts without central element $\kappa$ in (1.9) and (1.13), and introduce a notion of ``Lie algebra with one-variable structure.'' Let $V$ be a vector space and let $t$ again be an indeterminant. Form a tensor
$$\hat{V}=V\otimes_{\Bbb{F}}\Bbb{F}[t,t^{-1}].\eqno(1.14)$$
Define the {\it generating function}
$$u(z)=\sum_{n\in\Bbb{Z}}(u\otimes t^n)z^{-n-1}\qquad\for\;\;u\in V,\eqno(1.15)$$
where $z$ is again a formal variable. A {\it Lie algebra with one-variable structure} is a vector space $\hat{V}$ with the Lie bracket of the form
$$[u(z_1),v(z_2)]=\sum_{i=0}^m\sum_{j=0}^nz_2^{-1}\ptl_{z_1}^i\dlt\left({z_1\over z_2}\right)\ptl_{z_2}^jw_{ij}(z_2)\eqno(1.16)$$
for $u,v\in V$, where $m,\;n$ are nonnegative integers depending on $u,v$, and $w_{ij}\in V$.
Moreover, each $w_{ij}$ depends on $u$ and $v$ bilinearly by (1.15) and (1.16). 

 According to Wightman's axioms of quantum field theory,
the algebraic content of two-dimensional quantum field theory is a certain new representation theory (not the same as the traditional theory!) of Lie algebras with one-variable structure in more general sense (e.g., cf. [K1]), where intertwining operators among irreducible modules, ``partition functions'' (characters in algebraic terms), ``correlation functions'' related to (1.16) etc. play important roles. A ``conformal algebra'' is the local structure of the Lie algebras with one-variable structure. The above formula (1.16) can be represented by the following operator
$$Y^+(u,z)v=\sum_{i=0}^m\sum_{j=0}^n(-1)^ii!\ptl^jw_{ij}z^{-i-1},\eqno(1.17)$$
where ``$\ptl$'' is just a symbol standing for $\ptl_{z_2}$. 

\section{Basic Things about Conformal Algebras}

Before presenting the definition of conformal algebras, we shall introduce some notations. The following operator of taking residue will be used
$$\rd_z(\sum_{j\in\Bbb{Z}}\xi_jz^j)=\xi_{-1},\eqno(2.1)$$
where $\xi_j$ are in some vector space $V$. Moreover, all the binomials are assumed to be expanded in the second variable. For example,
$${1\over z-x}={1\over z(1-x/z)}=\sum_{j=0}^{\infty}z^{-1}\left({x\over z}\right)^j=\sum_{j=0}^{\infty}z^{-j-1}x^j.\eqno(2.2)$$
In particular, the above equation implies
$$\rd_x{1\over z-x}(\sum_{j\in\Bbb{Z}}\xi_jz^j)=\sum_{j=1}^{\infty}\xi_{-j}z^{-j}.\eqno(2.3)$$
So the operator $\rd_x(1/(z-x)$(---) is taking the part of negative powers in a formal series and changing the variable $x$ to $z$. For two vector spaces $V_1$ and $V_2$, we denote by $LM(V_1,V_2)$ the space of linear maps from $V_1$ to $V_2$. 
\psp 

 A {\it conformal algebra} $R$ is  an $\Bbb{F}[\ptl]$-module equipped with a  linear map $Y^+(\cdot,z):\;R\rightarrow LM(R,R[z^{-1}]z^{-1})$ satisfying:
$$Y^+(\ptl u,z)={d\over dz}Y^+(u,z),\eqno(2.4)$$
$$Y^+(u,z)v=\rd_x{1\over z-x}e^{x\ptl}Y^+(v,-x)u,\eqno(2.5)$$
$$Y^+(u,z_1)Y^+(v,z_2)-Y^+(v,z_2)Y^+(u,z_1)=\rd_x{1\over z_2-x}Y^+(Y^+(u,z_1-x)v,x) \eqno(2.6)$$
for $u,v\in R$. We denote by $(R,\ptl,Y^+(\cdot,z))$ a conformal algebra.
\psp

The above definition is the equivalent generating-function form to that given in [K1], where the author used the component formula with $Y^+(u,z)=\sum_{n=0}^{\infty}u_{(n)}z^{-n-1}/n!$. 
\psp

{\bf Remark 2.1}. (1) For $u\in R$,
$$Y(u,z)=\sum_{j=0}^{\infty}u_jz^{-j-1}\qquad\mbox{with}\;\;u_j\in\mbox{End}\:R,\eqno(2.7)$$
the algebra of linear transformations on $R$. Such an operator valued function has also been used in analysis. For any $u,v\in R$,
$$u_j(v)=0\qquad\for\;\;j>>0.\eqno(2.8)$$

(2) As an algebraic structure, conformal algebra had been essentially used by vertex operator algebraists for quite a few years before the name of conformal algebras was introduced by Kac [K1]. The notation $Y^+(\cdot,z)$ is used because a conformal algebra is essentially
``positive part'' of a vertex algebra. We refer to [X2] for the details of this point. 
 However, it is important to single out this structure. In terms of simple algebras, the category of conformal algebras is much smaller than that of vertex algebras. 

(3) A conformal algebra is equivalent to a linear Hamiltonian operator, a notion that introduced by Gel'fand and Dikii in mid seventies (see [X5] for the detailed proof).

(4) A notion of ``associative conformal algebra'' was introduced in [K2]. This notion seems the same as that of ``Gel'fand-Dorfman operator algebra'' in [X1], where I have proved that such an algebra with a certain unit element is equivalent to an associative algebra with a derivation.

(5) We refer to [X2] for twisted Lie algebras generated by conformal algebras and the detailed connections of conformal algebras with vertex operator algebras.

\psp

Note that a Lie algebra ${\cal G}$ can be viewed as a vector space with a linear map
$$\ad: {\cal G}\rta LM({\cal G},{\cal G})=\mbox{End}\:{\cal G}\eqno(2.9)$$
such that
$$\ad_u(v)=-\ad_v(u)\qquad(\mbox{Skew-Symmetry}),\eqno(2.10)$$
$$\ad_u\ad_v-\ad_v\ad_u=\ad_{\ad_{u}(v)}\qquad(\mbox{Jacobi Identity})\eqno(2.11)$$
for $u,v\in{\cal G}$. So the Expression (2.5) is an analogue of the skew-symmetry (2.10) and (2.6) is an analogue of the Jacobi identity (2.11).

We have the following obvious fact
$$\rd_z{d\over dz}(\sum_{j\in\Bbb{Z}}\xi_jz^j)=0.\eqno(2.12)$$
So
$$\rd_z{df(z)\over dz}g(z)=\rd_z{d(f(z)g(z))\over dz}-\rd_zf(z){dg(z)\over dz}=
-\rd_zf(z){dg(z)\over dz}.\eqno(2.13)$$
Using the above equation, I have given a direct elementary proof of the following theorem originally due to Kac (cf. [X5]).
\psp

{\bf Theorem 2.1 (Kac, [K1])}. {\it Let} $R$ {\it be an} $\Bbb{F}[\ptl]$-{\it module and let} $V$ {\it be a subspace of} $R$ {\it such that}
$$R=\Bbb{F}[\ptl]V.\eqno(2.14)$$
{\it Suppose that we have a linear map} $Y^+(\cdot,z): R\rta LM(R,R[z^{-1}]z^{-1})$ {satisfying}
$$[\ptl,Y^+(u,z)]=Y^+(\ptl u,z)={d\over dz}Y^+(u,z)\eqno(2.15)$$
{\it for} $u\in R$. {\it Then the triple} $(R,\ptl,Y^+(\cdot,z))$ {\it forms a conformal algebras if and only if Axiom (2.5), and Axiom (2.6) acting on} $V$ {\it hold for} $u,v\in V$. 
\psp

In fact, we have also used (2.13) to prove the converse theorem.
\psp

{\bf Theorem 2.2 (Xu, [X5])}. {\it For a conformal algebra} $(R,\ptl,Y^+(\cdot,z))$, {\it we always have}
$$[\ptl,Y^+(u,z)]=Y^+(\ptl u,z)\qquad\mbox{\it for}\;\;u\in R.\eqno(2.16)$$

Theorem 2.1 is very important in constructing conformal algebras.
\psp

{\bf Example 2.1}. Let ${\cal G}$ be a Lie algebra. Set
$$R({\cal G})=\Bbb{F}[\ptl]\otimes_{\Bbb{F}}{\cal G},\eqno(2.17)$$
which is a free $\Bbb{F}[\ptl]$-module. 
We identify $1\otimes{\cal G}$ with ${\cal G}$. Define
$$Y^+(u,z)v=[u,v]z^{-1}\qquad\for\;\;u,v\in{\cal G}.\eqno(2.18)$$
Expression (2.15) force us to extend $Y^+(\cdot,z)$ on $R$ uniquely by
$$Y^+(\ptl^mu,z)\ptl^nv=\sum_{j=0}^n(-1)^j\left(\!\!\begin{array}{c} n\\ j \end{array}\!\!\right)\left({d\over dz}\right)^{m+j}
\ptl^{n-j}Y^+(u,z)v\eqno(2.19)$$
for $m,n\in\Bbb{Z}$ and $u,v\in V$. So (2.4)  naturally holds. Axiom (2.5) is equivalent to the Lie algebra skew-symmetry and Axiom (2.6) is equivalent to the Jacobi identity. Thus $R({\cal G})$ forms a conformal algebra, which generates the loop Lie algebra ${\cal G}\otimes \Bbb{F}[t,t^{-1}]$ (cf. (1.16), (1.17)) with the Lie bracket
$$[u\otimes t^j,v\otimes t^l]=[u,v]\otimes t^{j+l}\qquad\for\;\;u,v\in{\cal G},\;j,l\in\Bbb{Z}.\eqno(2.20)$$

{\bf Example 2.2}. Let
$$R_W=\Bbb{F}[\ptl]e.\eqno(2.21)$$
Then we have the structure map determined by
$$Y^+(e,z)e=\ptl ez^{-1}+2ez^{-2}.\eqno(2.22)$$
Conformal algebra $R_W$ generates the rank-one Witt algebra ${\cal W}_1=\Bbb{F}e\otimes_\Bbb{F}\Bbb{F}[t,t^{-1}]$ (cf. (1.16), (1.17)) with the Lie bracket
$$[e\otimes t^j,e\otimes t^l]=(j-l)e\otimes t^{j+l}\qquad\for\;\;j,l\in\Bbb{Z}.\eqno(2.23)$$

An {\it ideal} ${\cal I}$ of a conformal algebra $R$ is a subspace such that
$$\ptl({\cal I})\subset{\cal I},\;\;Y^+(u,z)({\cal I})\subset {\cal I}[z^{-1}]z^{-1}.\eqno(2.24)$$
A conformal algebra $R$ is called of {\it finite type} if $R$ has a finite-dimensional subspace $V$ such that $R=\Bbb{F}[\ptl]V$.
\psp

{\bf Theorem 2.3 (D'Andrea and Kac, [DK])}. {\it A simple conformal algebra of finite type is either isomorphic to} $R_W$ {\it or isomorphic to} $R({\cal G})$ {\it for some finite-dimensional simple Lie algebra} ${\cal G}$.
\psp

The above theorem shows that there are no new algebras in simple conformal algebras of finite type. Zel'manov has an approach to simple conformal algebras of infinite type by introducing certain filtrations and Gel'fand-Kirillov dimensions. I have adopted the following two approaches to simple conformal algebras of infinite type:
\pse

1. Let $R$ be a conformal algebra that is a free $\Bbb{F}[\ptl]$-modules over a subspace $V$, namely,
$$R=\Bbb{F}[\ptl]V\cong \Bbb{F}[\ptl]\otimes_{\Bbb{F}} V.\eqno(2.25)$$
To my best knowledge, all the known simple conformal algebras have this property. Let $m$ be a positive integer. The algebra $R$ is said of {\it degree} $m$ if for any $u,v\in V$,
$$Y^+(u,z)v=\sum_{0<j;i+j\leq m}\ptl^iw_{i,j}z^{-j}\qquad\mbox{with}\;\;w_{i,j}\in V,\eqno(2.26)$$
and $w_{m-j,j}\neq 0$ for some $u,v\in V$ and $j\in\{1,2,...,m\}$.

A linear conformal algebra is isomorphic to $R({\cal G})$ in Example 2.1.
A quadratic conformal algebra is equivalent to a bialgebraic structure on $V$, which I called a {\it Gel'fand-Dorfman bialgebra} in [X3].
One of the structures is a Lie algebra, the other is a Novikov algebra and they satisfy a compatibility equation of five nonzero terms (cf. [X3]). We refer to
[X6] for a survey on Gel'fand-Dorfman bialgebras, and special families of simple cubic conformal algebras and simple quartic conformal algebras. 

This approach is to study simple conformal algebras of finite degree. The development of this approach is at initial stage.

\pse

2. A conformal algebra $R$ is called $\Bbb{N}$-{\it weighted} if
$$R=\bigoplus_{n=0}^{\infty}R^{(n)}\eqno(2.27)$$
is a graded space such that for $u\in R^{(l)}$, we have
$$Y^+(u,z)=\sum_{j=1-l}^{\infty}u(j)z^{-j-l},\;\;u(j)(R^{(n)})\subset R^{(n-j)}\eqno(2.28)$$
(cf. [X4]). Let $R$ be an $\Bbb{N}$-{\it weighted} conformal algebra that is a free $\Bbb{F}[\ptl]$-module over an $\Bbb{N}$-graded subspace $V$. The algebra
$R$ is called a conformal algebra of {\it finite growth} if
$$\dim V\bigcap R^{(n)}\leq N_0\qquad\for\;\;n\in\Bbb{N},\eqno(2.29)$$
where $N_0$ is a fixed positive integer. 

This second approach is to study simple conformal algebras of finite growth. Up to this point, we only know several families of such algebras (cf. [X4]). In next two sections, we shall introduce three families of simple conformal algebras of finite growth related to Jordan algebras, where the first family is directly taken from [X4], and the second and third are extended versions of those in [X4],

\section{Free Quadratic Fermionic Fields with Derivatives}

In this section, we shall present the conformal subalgebras related to quadratic free fermionic fields with derivatives.

Let $H$ be a vector space with a nondegenerate symmetric bilinear form $\la\cdot,\cdot\ra$ such that there exist two subspaces $H_+,H_-$ satisfying $H=H_++H_-$ and
$$\la H_+,H_+\ra=\la H_-,H_-\ra=\{0\}.\eqno(3.1)$$
Thus $H_+$ is isomorphic to the dual space $(H_-)^{\ast}$ through the nondegenerate symmetric bilinear form $\la\cdot,\cdot\ra$. 

Let $t$ be an indeterminate and set
$$\hat{H}=H\otimes_{\Bbb{F}}\Bbb{F}[t,t^{-1}]t^{1/2}\oplus \Bbb{F}\kappa,\eqno(3.2)$$
where $\kappa$ is a symbol to denote a base vector of one-dimensional vector space. We define algebraic operation $[\cdot,\cdot]$ on $\hat{H}$ by
$$[h_1\otimes t^m+\lmd_1\kappa,h_2\otimes t^n+\lmd_2\kappa]=\la h_1,h_2\ra\dlt_{m+n,0}\kappa\eqno(3.3)$$
for $h_1,h_2\in H,\;m,n\in\Bbb{Z}+1/2,\;\lmd_1,\lmd_2\in\Bbb{F}$. Then $(\hat{H},[\cdot,\cdot])$ forms a Lie superalgebra with the $\Bbb{Z}_2$-grading
$$\hat{H}_1=H\otimes_{\Bbb{F}}\Bbb{F}[t,t^{-1}]t^{1/2},\;\;\hat{H}_0=\Bbb{F}\kappa.\eqno(3.4)$$
 For convenience, we denote
$$h(m)=h\otimes t^m\qquad\for\;\;h\in H,\;m\in\Bbb{Z}+1/2.\eqno(3.5)$$
Set
$$\hat{H}_-=\mbox{span}\:\{h(- m)\mid h\in H,\;m\in\Bbb{N}+1/2\},\eqno(3.6)$$
$$\hat{B}_H=\mbox{span}\:\{\kappa,h(m)\mid h\in H,\;m\in\Bbb{N}+1/2\}.\eqno(3.7)$$
Then $\hat{H}_-$ and $\hat{B}_H$ are trivial Lie sub-superalgebras of $\hat{H}$ and
$$\hat{H}=\hat{H}_-\oplus \hat{B}_H.\eqno(3.8)$$

Let $\Bbb{F}{\bf 1}$ be a one-dimensional vector space with the base element ${\bf 1}$. We define an action of $\hat{B}_H$ on 
$\Bbb{F}{\bf 1}$ by
$$h(m)({\bf 1})=0,\;\;\kappa({\bf 1})={\bf 1}\qquad\for\;\;h\in H,\;m\in\Bbb{N}+1/2.\eqno(3.9)$$
Then $\Bbb{F}{\bf 1}$ forms a $\hat{B}_H$-module. We denote by $U(\cdot)$ the universal envelopping algebra of a Lie algebra and by $\bigwedge(\cdot)$ the exterior algebra generated by a vector space. Form an induced $\hat{H}$-module
$${\cal F}=U(\hat{H})\otimes_{U(\hat{B}_H)}\Bbb{F}{\bf 1}\;\;(\cong \bigwedge(\hat{H}_-)\otimes_{\Bbb{F}}\Bbb{F}{\bf 1}\;\mbox{as vector spaces}).\eqno(3.10)$$
Moreover, we set
$$h^+(z)=\sum_{m=0}^{\infty}h(m+1/2)z^{-m-1},\;\;\;h^-(z)=\sum_{m=1}^{\infty}h(-m+1/2)z^{m-1},\eqno(3.11)$$
$$h(z)=h^+(z)+h^-(z)\eqno(3.12)$$
for $h\in H$. As operators on ${\cal F}$, $\{h(z)\mid h\in H\}$ are called {\it free fermionic fields}.

For convenience, we denote
$$u\otimes{\bf 1}=u\qquad\for\;\;u\in \bigwedge(\hat{H}_-).\eqno(3.13)$$
Set
$$\hat{R}_2=\mbox{span}\:\{h_1(-m)h_2(-n),{\bf 1}\mid h_1\in H_+,\;h_2\in H_-,\;m,n\in\Bbb{N}+1/2\}.\eqno(3.14)$$
We define a linear map $Y(\cdot,z):\hat{R}_2\rightarrow LM({\cal F},{\cal F}[z^{-1},z]])$ by
$$Y(h_1(-m-1/2)h_2(-n-1/2),z)={1\over m!n!}\left({d^mh^-_1(z)\over dz^m}{d^nh_2(z)\over dz^n}-{d^nh_2(z)\over dz^n}{d^mh^+_1(z)\over dz^m}\right)\eqno(3.15)$$
for $h_1\in H_+,\;h_2\in H_-$ and $m,n\in\Bbb{N}$ and
$$Y({\bf 1},z)=\mbox{Id}_{\cal F}.\eqno(3.16)$$
The operator $Y(h_1(-m-1/2)h_2(-n-1/2),z)$ is  a {\it quadratic fermionic field with derivatives}. Moreover, we write
$$Y(u,z)=\sum_{n\in\Bbb{Z}}u_nz^{-n-1},\;\;Y^+(u,z)=\sum_{n=0}^{\infty}u_nz^{-n-1}\qquad\for\;\;u\in\hat{R}_2.\eqno(3.17)$$
In particular, 
\begin{eqnarray*}& &(h_1(-m-1/2)h_2(-n-1/2))_k\\&=&\sum_{j=0}^{\infty}[\left(\!\!\begin{array}{c}-j-1\\ n\end{array}\!\!\right)\left(\!\!\begin{array}{c}j+m+n-k\\ m\end{array}\!\!\right)h_1(k-m-n-j-1/2)h_2(j+1/2)\\& &-\left(\!\!\begin{array}{c}-j-1 \\ m\end{array}\!\!\right)\left(\!\!\begin{array}{c}j+m+n-k\\ n\end{array}\!\!\right)h_2(k-m-n-j-1/2)h_1(j+1/2)]\hspace{1.5cm}(3.18)\end{eqnarray*}
for $h_1\in H_+,\;h_2\in H_-$ and $m,n,k\in\Bbb{N}$.

 Note that
$$h_1(m)h_2(n)(h_3(-j)h_4(-k))=\dlt_{m,k}\dlt_{n,j}\la h_1,h_4\ra \la h_2,h_3\ra{\bf 1},\eqno(3.19)$$
$$h_1(-m)h_2(n)(h_3(-j)h_4(-k))=\dlt_{n,j}\la h_2,h_3\ra h_1(-m)h_4(-k),\eqno(3.20)$$
$$h_2(-m)h_1(n)(h_3(-j)h_4(-k))=\dlt_{n,k}\la h_1,h_4\ra h_3(-j)h_2(-m)\eqno(3.22)$$
for $h_1,h_3\in H_+,\;h_2,h_4\in H_-$ and $j,k,m,n\in\Bbb{N}+1/2$. Expressions (3.18)-(3.22) show that
$$Y^+(u,z)v\subset \hat{R}_2[z^{-1}]z^{-1}\qquad\for\;\;u,v\in \hat{R}_2.\eqno(3.23)$$
Moreover, we define $\ptl\in\Edo \hat{R}_2$ by $\ptl({\bf 1})=0$ and
\begin{eqnarray*}\ptl (h_1(-m-1/2)h_2(-n-1/2))&=&(m+1)h_1(-m-3/2)h_2(-n-1/2)\\& &+(n+1)h_1(-m-1/2)h_2(-n-3/2)\hspace{1.8cm}(3.24)\end{eqnarray*}
for $h_1,h_2\in H$ and $m,n\in\Bbb{N}$. Then the family $(\hat{R}_2,\ptl,Y^+(|_{\hat{R}_2},z))$ forms a conformal algebra.

According to linear algebra, there exist basis $\{\vs^{\pm}_j|j\in I\}$ of $H_{\pm}$ (cf. (3.1)) such that
$$\la\vs^+_i,\vs_j^-\ra=\dlt_{i,j}\qquad\for\;\;i,j\in I\eqno(3.25)$$
by (3.1) and nondegeneracy of $\la\cdot,\cdot\ra$, where $I$ is an index set. Note that
$$\vs^+_{j_1}(-m)\vs^-_{j_2}(m)(\vs^+_{j_3}(-m)\vs^-_{j_4}(-m))=\dlt_{j_2,j_3}\vs^+_{j_1}(-m)\vs^-_{j_4}(-m)\eqno(3.26)$$
for $j_1,j_1,j_3,j_4\in I$ and $m\in\Bbb{N}+1/2$. Expression (3.26) is essentially equivalent to matrix multiplications! This gives us a motivation of constructing simple conformal algebras in next section.

\section{Conformal Algebras Generated by Jordan Algebras}

Let $k$ be a positive integer and set
$${\cal A}=M_{k\times k} (\Bbb{F}),\eqno(4.1)$$
the algebra of $k\times k$ matrices. Form a tensor
$$R_{k\times k}={\cal A}\otimes_{\Bbb{F}}\Bbb{F}[t_1,t_2],\eqno(4.2)$$
where $t_1,t_2$ are indeterminants. For convenience, we denote
$$u(m,n)=u\otimes t_1^mt_2^n\qquad\for\;\;u\in{\cal A},\;m,n\in\Bbb{N}.\eqno(4.3)$$
We define the $\Bbb{F}[\ptl]$-action on $R_{k\times k}$ by
$$\ptl(u(m,n))=(m+1)u(m+1,n)+(n+1)u(m,n+1)\eqno(4.4)$$
and the structure map $Y^+(\cdot,z)$ on $R_{k\times k}$ by
\begin{eqnarray*}& &Y^+(u(m_1,m_2),z)v(n_1,n_2)\\&=&
\left(\!\!\begin{array}{c} -n_1-1\\ m_2\end{array}\!\!\right)
\sum_{p=0}^{m_1+m_2+n_1}\left(\!\!\begin{array}{c}p\\ m_1 \end{array}\!\!\right)(uv)(p,n_2)z^{p-m_1-m_2-n_1-1}\\&&-\left(\!\!\begin{array}{c}-n_2-1\\ m_1\end{array}\!\!\right)\sum_{q=0}^{m_1+m_2+n_2}\left(\!\!\begin{array}{c} q\\ m_2\end{array}\!\!\right)(vu)(n_1,q)z^{q-m_1-m_2-n_2-1}\hspace{4cm}(4.5)\end{eqnarray*}
for $u,v\in{\cal A}$ and $m,m_1,m_2,n,n_1,n_2\in\Bbb{N}$. The above formula was motivated by the following formula 
of the nonnegative operators of a quadratic field acting on a quadratic element:
$$Y^+(\vs_{j_1}^+(-m_1-1/2)\vs_{j_2}^-(-m_2-1/2),z)\vs_{j_3}^+(-n_1-1/2)\vs^-_{j_4}(-n_2-1/2) \eqno(4.6)$$
for $j_1,j_2,j_3,j_4\in I$ and $m_1,n_1,m_2,n_2\in\Bbb{N}$, which is defined by (3.1), (3.15), (3.17) and (3.25). The reader may calculate (4.6) by (3.18) and compare it with (4.5) in order to better understand (4.5).

 We define weights by
$$R_{k\times k}^{(n)}=\mbox{span}\:\{u(m_1,m_2)\mid u\in{\cal A},\;m_1,m_2\in\Bbb{N},\;m_1+m_2+1=n\}\eqno(4.7)$$
for $n\in\Bbb{N}+1$. Then $(R_{k\times k},\ptl,Y^+(\cdot,z))$ forms a simple
conformal algebra of finite growth (cf. [X4])  with
$$V={\cal A}\otimes_{\Bbb{F}}\Bbb{F}[t_2]=\mbox{Span}\:\{u(0,n)\mid u\in{\cal A},\;n\in\Bbb{N}\}\eqno(4.8)$$
(that is, $R_{k\times k}$ is an $\Bbb{F}[\ptl]$-module over $V$). Next we want to pick out simple subalgebras.

For $\ell\in\Bbb{N}$, we set
$$R_{k\times k,\ell+1}=\mbox{Span}\:\{u(m,n+\ell)\mid u\in{\cal A},\;m,n\in\Bbb{N}\},\eqno(4.9)$$
\begin{eqnarray*}R_{k\times k,\ell+1}^{\ast}&=&\mbox{Span}\:\{\left(\!\!\begin{array}{c}n+\ell \\ \ell\end{array}\!\!\right)u(m,n+\ell)-(-1)^{\ell}\left(\!\!\begin{array}{c} m+\ell\\ \ell\end{array}\!\!\right)u^T(n,m+\ell)\\& &\mid u\in{\cal A},\;m,n\in\Bbb{N}\},\hspace{8.9cm}(4.10)\end{eqnarray*}
where the up index ``T'' stands for the transpose of matrices. When $k=2k_1$ is an even integer, we have another involutive anti-automorphism ``$\dg$'' of ${\cal A}$:
$$ u^{\dg}=\left(\begin{array}{rr}&-I_{k_1}\\ I_{k_1}&\end{array}\right)u^T\left(\begin{array}{rr}&I_{k_1}\\ -I_{k_1}&\end{array}\right)\qquad\for\;\;u\in{\cal A},\eqno(4.11)$$
where $I_{k_1}$ is the $k_1\times k_1$ identity matrix. Furthermore, we set
\begin{eqnarray*}R_{k\times k,\ell+1}^{\dg}&=&\mbox{Span}\:\{\left(\!\!\begin{array}{c}n+\ell \\ \ell\end{array}\!\!\right)u(m,n+\ell)-(-1)^{\ell}\left(\!\!\begin{array}{c} m+\ell\\ \ell\end{array}\!\!\right)u^{\dg}(n,m+\ell)\\& &\mid u\in{\cal A},\;m,n\in\Bbb{N}\}.\hspace{8.9cm}(4.12)\end{eqnarray*}

{\bf Theorem 4.1 (Xu, [X4])}. {\it For any} $1\leq \ell\in\Bbb{N}$, {\it all the subspaces} $R_{k\times k,\ell},\;R_{k\times k,\ell}^{\ast}$ {\it and} $R_{k\times k,\ell}^{\dg}$ {\it form simple conformal subalgebras of finite growth with the minimal weight} $\ell$.
\psp

A subset $S$ of a conformal algebra $R$ is called a {\it generator subset} of $R$ if
$$R=\mbox{Span}\:\{u^1_{n_1}u^2_{n_2}\cdots u_{n_s}^s(v)\mid u^i,v\in S,\;s,n_j\in\Bbb{N}\},\eqno(4.13)$$
where we have used the notation
$$Y^+(u,z)=\sum_{n=0}^{\infty}u_nz^{-n-1}\qquad\for\;\;u\in R.\eqno(4.14)$$

{\bf Theorem 4.2 (Xu, [X4])}. {\it For any} $2\leq \ell\in\Bbb{N}$, {\it all the conformal algebras} $R_{k\times k,\ell}\;(k\geq 2),\;R_{k\times k,\ell}^{\ast}\;(k\geq 2)$ {\it and} $R_{k\times k,\ell}^{\dg}\;(k\geq 4)$ {\it are generated by their subspace of minimal weight}.
\psp

Let
$$X_{\ell}=R_{k\times k,\ell},\;R_{k\times k,\ell}^{\ast}\;\;\mbox{or}\;\;R_{k\times k,\ell}^{\dg}.\eqno(4.15)$$
For any
$$u\in X_{\ell}^{(\ell)}=X_{\ell}\bigcap R_{k\times k}^{(\ell)},\eqno(4.16)$$
we write
$$Y^+(u,z)=\sum_{n=1-\ell}^{\infty}u(n)z^{-n-\ell}.\eqno(4.17)$$
Then $u(n)\in \mbox{End}\:X_{\ell}$ is an operator of degree $-n$. In particular, we can define an algebraic operation $\circ$ on its subspace $X_{\ell}^{(\ell)}$ of minimal weight by 
$$u\circ v=u(0)v\qquad\for\;\;u,v\in X_{\ell}^{(\ell)}.\eqno(4.18)$$

{\bf Theorem 4.3 (Xu, [X4])}. {\it When} $\ell$ {\it is even, the algebra} $(X^{(\ell)}_{\ell},\circ)$ {\it forms a simple Jordan algebra of type A, B or C, which is independent of} $\ell$.
\psp

{\bf Remark 4.4}. (1) When $\ell$  is odd, the algebra $(X^{(\ell)}_{\ell},\circ)$  forms a Lie algebra isomorphic to $gl(k),\;o(k)$ or $sp(k)$,  which is independent of $\ell$.

(2) The algebras $R_{k\times k,\ell}^{\ast}$ and $R_{k\times k,\ell}^{\dg}$ with $\ell>2$ have not been defined in [X4]. However, the proofs of Theorems 3.1, 3.2 and 3.3 in [X4] show that the same conclusions hold for them.

(3) The other cases such as $\ell=1$ in Theorem 4.2 have been handled in [X4].
\vspace{1cm}

\noindent{\Large \bf References}

\hspace{0.5cm}

\begin{description}

\item[{[DK]}] A. D'Andrea and V. G. Kac, Structure theory of finite conformal algebras, {\it Selecta Math.} (N.S.) {\bf 4} (1998), 377-418.

\item[{[K1]}] V. G. Kac, {\it Vertex algebras for beginners}, University lectures series, Vol {\bf 10}, AMS. Providence RI, 1996.

\item[{[K2]}] ---, Formal distribution algebras and conformal algebras, {\it XIIth International Congress of Mathematical Physics (ICMP '97) (Brisbane), 80-97, Internat. Press, Cambridge, MA, 1999}.

\item[{[LW]}] J. Lepowsky and R. L. Wilson, Construction of affine Lie algebra
$A^{(1)}$, {\it Commun. Math. Phys.} {\bf 62} (1978), 43-53.

\item[{[X1]}] X. Xu, Hamiltonian operators and associative algebras with a derivation, {\it Lett. Math. Phys.} {\bf 33} (1995), 1-6.

\item[{[X2]}] ---, {\it Introduction to vertex operator superalgebras and their modules}, Kluwer Academic Publishers, Dordrecht/Boston/London, 1998.

\item[{[X3]}] ---,  Quadratic conformal superalgebras, {\it J. Algebra} {\bf 231} (2000), 1-38.

\item[{[X4]}] ---, Simple conformal superalgebras of finite growth, {\it Algebra Collquium} {\bf 7} (2000), 205-240.

\item[{[X5]}]---, Equivalence of conformal superalgebras to Hamiltonian superoperators,  to appear in {\it Algebra Colloquium} {\bf 8} (2001), 63-92.

\item[{[X6]}]---, Gel'fand-Dorfman Bialgebras, {\it arXiv:math.QA/0008223}.

\end{description}

\noindent Department of Mathematics, The Hong Kong University of Science and Technology, Clear Water Bay, Kowloon, Hong Kong

\end{document}